\documentclass[twoside,11pt]{article}

\usepackage{blindtext}

\usepackage[preprint]{jmlr2e}

\usepackage[OT2, T1]{fontenc}
\usepackage[russian, english]{babel}
\usepackage{amsmath}
\usepackage{hyperref}
\hypersetup{colorlinks=true,linkcolor=blue,citecolor=blue}

\newcommand{\load}{.}
\usepackage{\load/mystyle}

\usepackage{lastpage}
\jmlrheading{23}{2022}{1-\pageref{LastPage}}{1/21; Revised 5/22}{9/22}{21-0000}{Guillaume Garrigos}

\ShortHeadings{Square distance functions are Polyak-Łojasiewicz and vice-versa}{G. Garrigos}
\firstpageno{1}

\begin{document}

\title{Square distance functions are Polyak-Łojasiewicz\\ and vice-versa}

\author{\name Guillaume Garrigos \email garrigos@lpsm.paris \\
       \addr Université Paris Cité and Sorbonne Université, CNRS\\
       Laboratoire de Probabilités, Statistique et Modélisation\\
       F-75013 Paris, France}

\editor{My editor}

\maketitle

\begin{abstract}%   <- trailing '%' for backward compatibility of .sty file
This short note gathers known results to state that the squared distance function to a (nonconvex) closed set of an Euclidean space is Polyak-Łojasiewicz.
As a fuzzy reciprocate, we also recall that every Polyak-Łojasiewicz function can be bounded from below by the squared distance function to its set of minimizers.
\end{abstract}


\section{Introduction}

The Polyak-Łojasiewicz (PŁ) property (see Definition \ref{D:PL} below) appeared in a work by Polyak \citeyearpar[\foreignlanguage{russian}{Teorema 4}]{Pol63}, in parallel to the development of Łojasiewicz's inequalities \citeyearpar[Théorème 4]{Loj63a} for the study of analytic gradient flows.
The machine learning optimization community recently gained interest for this property \citep{KarNutSch16}, because making this hypothesis allows to obtain results in a nonconvex setting which are quantitatively similar to what is usually obtained in a strongly convex setting \citep{FosSekSri18,LeiHuLiTan19}.

While it is generally emphasized that the PŁ property can be derived from strong convexity (or other weakened notions), most works seem to elude the question of what are \emph{nonconvex} Polyak-Łojasiewicz functions, except for a few toy examples.
Note that we are considering the original PŁ inequality, which is \emph{global} by nature, and not ``local'' PŁ inequalities which can be derived for some machine learning models \citep{HarMa16,LiYua17,LiuZhuBel22}.
%Those are more delicate to exploit since they require the initialization of optimization algorithm to be close enough from a solution.

In this letter we put the emphasis on a class of functions which we show to be PŁ: the squared distance to closed sets (Theorem \ref{T:distance squared PL}).
Not only this provides a class of functions illustrating the PŁ property, which can have an interesting pedagogical purpose, but we also show that every PŁ function can be bounded from below (and above, if smooth) by a squared distance function (Corollary~\ref{T:PL implies 2 conditioning}).
The first result does not seem to appear in the literature, even though it is a simple application of standard results in variational analysis.
The second was essentially proved by \cite{BolDanLeyMaz10}, but might have been missed by the machine learning community, and its proof is adapted here to fit our needs.

\section{Definitions and notations}

Given a closed set $\Omega \subset \mathbb{R}^N$, we note $d_\Omega : \mathbb{R}^N \to \mathbb{R}$ the \emph{distance function to $\Omega$}
\begin{equation*}
(\forall x \in \mathbb{R}^N) \quad
d_\Omega(x) := \inf\limits_{p \in \Omega} \Vert p - x \Vert,
\end{equation*}
and $\proj_\Omega : \mathbb{R}^N 
\rightrightarrows
\mathbb{R}^N$ the \emph{projection mapping onto $\Omega$}
\begin{equation*}
(\forall x \in \mathbb{R}^N) \quad
\proj_\Omega(x) := \underset{p \in \Omega}{\rm{argmin}}~\Vert p - x \Vert.
\end{equation*}
Note that we do not assume $\Omega$ to be convex, so the projection $\proj_\Omega(x)$ could be a set larger than a singleton. But it is always nonempty, since the function $p \mapsto \Vert p-x \Vert$ is coercive and continuous.
We recall if needed that $d_\Omega$ is Lipschitz continuous \cite[Example 9.6]{RocWet09}.

Given a function $f : \mathbb{R}^N \to \mathbb{R}$, and $x \in \mathbb{R}^N$, we note $\partialF f(x)$  the \emph{Fréchet subdifferential of $f$ at $x$}  defined by
\begin{equation*}
x^* \in \partialF  f(x) \Leftrightarrow \liminf\limits_{y \to x} \frac{f(y) - f(x) - \langle x^*, y-x \rangle}{\Vert y - x \Vert} \geq 0.
\end{equation*}
Moreover, we note $\partialL f(x)$ the \emph{Limiting subdifferential of $f$ at $x$}, which is composed of vectors $x^*\in \partialL  f(x)$ for which
%such that there exists sequences $x_n \to x$ and $x_n^* \to x^*$ such that $x_n^* \in \partialF  f(x_n)$ and $f(x_n) \to f(x)$.
\begin{equation*}
(\exists x_n \to x)(\exists x_n^* \to x^*) \quad f(x_n) \to f(x) \text{ and } x_n^* \in \partialF  f(x_n).
\end{equation*}
It is immediate from the  definitions that $\partialL  (\alpha f)(x) = \alpha \partialL  f(x)$ for any $\alpha \geq 0$ and that $\partialF  f(x) \subset \partialL  f(x)$.
We recall if needed that if $f$ is continuously differentiable at $x$, then $\partialF  f(x) = \partialL  f(x) = \{\nabla f(x) \}$ \cite[Exercise 8.8]{RocWet09}; and that if $f$ is convex then $\partialF f(x)$ and $\partialL f(x)$ both coincide with the classical convex subdifferential \cite[Proposition 8.12]{RocWet09}.
Finally, we will say that $x$ is a \emph{Limiting critical point} of $f$ if $0 \in \partialL f(x)$.

\bigskip

We now turn to the main property of interest in this document.
Because square distance functions to nonconvex sets are nonsmooth, we have to consider a nonsmooth version of the PŁ property.
We are implicitly making choices when doing so, which we  justify in Remark~\ref{R:choice of subdifferential}.

\begin{definition}\label{D:PL}
Let $f : \mathbb{R}^N \to \mathbb{R}$.
We say that $f$ is \textbf{Polyak-Łojasiewicz} (\textbf{PŁ} for short) if there exists $\mu>0$ such that
\begin{equation*}
(\forall x \in \mathbb{R}^N)(\forall x^* \in \partialL  f(x)) \quad f(x) - \inf f \leq \frac{1}{2\mu} \Vert x^* \Vert^2.
\end{equation*}
%We note $Ł^2_\mu(\mathbb{R}^N)$ the space of PŁ functions with constant $\mu$.
\end{definition}

\begin{figure}
\label{F:distance parabola}
\begin{center}
\includegraphics[width=0.5\linewidth]{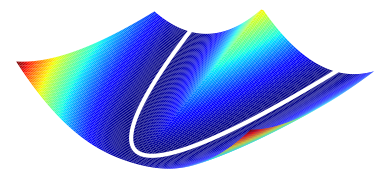}
\end{center}
\caption{Graph of $f(x) = \frac{1}{2}d_\Omega(x)^2$, where $\Omega = \{(t,t^2)\in \mathbb{R}^2 \ | \ t \in \mathbb{R}\}$.}
\end{figure}

\begin{remark}[Invexity]
The PŁ property is not related to convexity (see for instance Figure~\ref{F:distance parabola} for a nonconvex PŁ function), but shares a strong global property with convex functions, which is that every Limiting critical point must be a global minimizer.
See Example \ref{Ex:PL clarke subgradients} on why the word \emph{Limiting} is important here.
\end{remark}

\noindent The PŁ property describes functions having a global quadratic growth with respect to its minimizer set.
One can also consider functions having a more general polynomial growth:

\begin{definition}\label{D:Łojasiewicz global}
Let $f : \mathbb{R}^N \to \mathbb{R}$, $p>1$ and $\mu >0$.
We say that $f$ is \textbf{globally $p$-Łojasiewicz} if there exists $\mu>0$ such that
\begin{equation*}
(\forall x \in \mathbb{R}^N)(\forall x^* \in \partialL  f(x)) \quad f(x) - \inf f \leq \frac{1}{q\mu^{\frac{q}{p}}} \Vert x^* \Vert^q,
\quad \text{ where } \frac{1}{p}+ \frac{1}{q} = 1.
\end{equation*}
%We note $Ł^p_\mu(\mathbb{R}^N)$ the space of globally $p$-Łojasiewicz with constant $\mu$.
\end{definition}

\begin{example}[Power norms]
Let $f(x) = \frac{\mu}{p} \Vert x \Vert^p$ for some $p>1$ and $\mu>0$.
%Then it is easy to verify that $f \in Ł^p_\mu(\mathbb{R}^N)$, i.e. it is globally $p$-Łojasiewicz.
Then it is easy to verify that $f$ is globally $p$-Łojasiewicz with constant $\mu$.
\end{example}

\begin{example}[Uniformly convex functions]
Let $f : \mathbb{R}^N \to \mathbb{R}$ be $p$-uniformly convex for $p \geq 2$, in the sense that there exists $\mu>0$ such that for all $x,y \in \mathbb{R}^N$ and all $\alpha \in [0,1]$,
\begin{equation*}
f((1-\alpha)x + \alpha y) + \alpha(1-\alpha) \frac{\mu}{p} \Vert y-x \Vert^p \leq (1-\alpha)f(x) + \alpha f(y).
\end{equation*}
%Then $f \in Ł^p_\mu(\mathbb{R}^N)$, i.e. it is globally $p$-Łojasiewicz (see e.g. \cite[Lemma A.6]{GarRosVil22}).
Then $f$ is globally $p$-Łojasiewicz  with constant $\mu$ \cite[Lemma A.6]{GarRosVil22}.
When $p=2$ we recover the well-known fact that strongly convex functions are Polyak-Łojasiewicz.
\end{example}

\section{The subdifferential of distance functions}

\begin{lemma}\label{L:subdifferential distance}
Let $\Omega \subset \mathbb{R}^N$ be closed and nonempty.
Then
\begin{enumerate}
	\item\label{L:subdifferential distance:inside} for every $x \in \Omega$, $0 \in \partialL  d_\Omega (x)$;
	\item\label{L:subdifferential distance:outside} for every $x \notin \Omega$, $\partialL  d_\Omega (x) = \frac{x - \proj_\Omega(x)}{d_\Omega(x)}$;
	\item\label{L:subdifferential distance:nonempty} for every $x \in \mathbb{R}^N$, $\partialL  d_\Omega(x) \neq \emptyset$.
\end{enumerate}
\end{lemma}

\begin{proof}
For item \ref{L:subdifferential distance:inside}, note that $x$ is a minimizer of $d_\Omega$, so the optimality condition $0 \in \partialL  d_\Omega(x)$ holds, see \cite[Theorem 10.1]{RocWet09}.
%
%take $x \in \Omega$ and use the fact that $d_\Omega(x) =0$ to see that with $x^*=0$ we have
%\begin{equation*}
%(\forall y \in \mathbb{R}^N \setminus \{x\}) \quad \frac{d_\Omega(y) - d_\Omega(x) - \langle x^*, y-x \rangle}{\Vert y-x \Vert} \geq 0.
%\end{equation*}
%This implies from the definition that $0 \in \partialF  d_\Omega (x) \subset \partialL  d_\Omega(x)$.
For item \ref{L:subdifferential distance:outside}, the proof can be found in \cite[Example 8.53]{RocWet09}.
Item \ref{L:subdifferential distance:nonempty} is a direct consequence of the two previous items, and the fact that the projection on a nonempty closed set is always nonempty.
\end{proof}

Using a suitable chain rule, we obtain a similar result for powers of distance function.

\begin{lemma}\label{L:subdifferential p distance}
Let $\Omega \subset \mathbb{R}^N$ be closed and nonempty.
Let $p > 1$, $\mu>0$ and let $f(x) := \frac{\mu}{p} d_\Omega(x)^p$.
Then
\begin{enumerate}
	\item for all $x \in \Omega$, $ \partialL  f(x) = \{0 \}$;
	\item for all $x \notin \Omega$, $\partialL  f(x) =\mu d_\Omega(x)^{p-2}({x - \proj_\Omega(x)}) $;
	\item for every $x \in \mathbb{R}^N$, $\partialL  f(x) \neq \emptyset$.
\end{enumerate}
\end{lemma}

\begin{proof}
Let us write $f = \phi \circ d_\Omega$, where $\phi(t) = \frac{\mu}{p}t^p$.
Using the chain rule\footnote{To help with the notations in \cite{Mor12}: one must look at \cite[Eq. 1.61]{Mor12} and take $\varphi(x,y) := \phi(y)$ which is strictly differentiable for $p>1$.} in \cite[Theorem 1.110.ii]{Mor12}, we can write that
\begin{equation*}
\partialL  f(x) = \partialL  \left( \phi'(d_\Omega(x)) d_\Omega(\cdot) \right)(x).
\end{equation*}
Here $\phi'(d_\Omega(x)) = \mu d_\Omega(x)^{p-1}$, so
\begin{equation*}
\partialL  f(x) = \partialL  \left( \mu d_\Omega(x)^{p-1} d_\Omega(\cdot) \right)(x)
=  \mu d_\Omega(x)^{p-1} \partialL  d_\Omega (x).
\end{equation*}
Now, use Lemma \ref{L:subdifferential distance} to see that $\partialL  d_\Omega (x) \neq \emptyset$.
So for $x \in \Omega$ we have  $d_\Omega(x)^{p-1} =0$ (because $p>1$) which implies that $ \partialL  f(x) = \{0 \}$.
Instead, if $x \neq \emptyset$, we can use Lemma \ref{L:subdifferential distance} and the fact that $d_\Omega(x) \neq 0$ to conclude that $ \partialL  f(x) = \mu d_\Omega(x)^{p-2}({x - \proj_\Omega(x)}) $.
\end{proof}

\section{Łojasiewicz inequalities for distance functions}

\begin{theorem}\label{T:distance function p loja}
Let $\Omega \subset \mathbb{R}^N$ be closed and nonempty.
Let $p > 1$, $\mu>0$ and let $f(x) := \frac{\mu}{p} d_\Omega(x)^p$.
%Then $f \in Ł^p_\mu(\mathbb{R}^N)$, i.e. it is globally $p$-Łojasiewicz.
Then $f$ is globally $p$-Łojasiewicz with constant $\mu$.
\end{theorem}

\begin{proof}
Let $x \in \mathbb{R}^N$ and let us verify that the Łojasiewicz inequality in Definition \ref{D:Łojasiewicz global} holds true.
Note that here $\inf f = 0$.
If $x \in \Omega$ then the inequality is trivial, because $f(x)=0$ and from Lemma \ref{L:subdifferential p distance} we know that $\partialL  f(x) = \{0\}$.
If $x \notin \Omega$, then we use Lemma \ref{L:subdifferential p distance} to write, for all $x^* \in \partialL  f(x)$,
\begin{equation*}
\Vert x^* \Vert^q = \mu^q d_\Omega(x)^{(p-2)q} d_\Omega(x)^q = \mu^q d_\Omega(x)^p = \mu^{q-1} p f(x),
\end{equation*}
from which we see that Definition \ref{D:Łojasiewicz global} is verified.
\end{proof}

\begin{theorem}\label{T:distance squared PL}
Let $\Omega \subset \mathbb{R}^N$ be closed and nonempty.
Let $\mu>0$ and let $f(x) := \frac{\mu}{2} d_\Omega(x)^2$.
%Then $f \in Ł^2_\mu(\mathbb{R}^N)$, i.e. it is Polyak-Łojasiewicz.
Then $f$ is Polyak-Łojasiewicz with constant $\mu$.
\end{theorem}

\begin{proof}
Apply Theorem \ref{T:distance function p loja} with $p=q=2$.
\end{proof}

Using similar arguments, we can also see that (powers of) distance functions verify two other related properties, which are typically used for describing the geometry of a convex function. 
See \citep{BolNguPeySut17,GarRosVil22} for more details and references on those properties in a convex setting.

\begin{proposition}\label{T:distance function p conditioning metric regular}
Let $\Omega \subset \mathbb{R}^N$ be closed and nonempty.
Let $p > 1$, $\mu>0$ and let $f(x) := \frac{\mu}{p} d_\Omega(x)^p$.
Then 
\begin{enumerate}
	\item $f$ is \textbf{globally $p$-conditioned} with constant $\mu$ :
	\begin{equation*}\label{eq:conditioning}
	(\forall x \in \mathbb{R}^N) \quad \frac{\mu}{p}d_{{\rm{argmin}} f}(x)^p \leq f(x) - \inf f;
	\end{equation*}
	\item $f$ is \textbf{globally $p$-submetric regular}  with constant $\mu$ :
	\begin{equation*}
	(\forall x \in \mathbb{R}^N)(\forall x^* \in \partialL  f(x)) \quad {\mu}d_{{\rm{argmin}} f}(x)^{p-1} \leq \Vert x^* \Vert.
	\end{equation*}
\end{enumerate}
\end{proposition}

\begin{proof}
The result follows again from Lemma \ref{L:subdifferential p distance}.
If $x \in \Omega$ the inequalities are trivial, and if $x \notin \Omega$ we use the fact that
$d_{{\rm{argmin}} f}(x) = d_\Omega(x)$, $f(x) - \inf f = \frac{\mu}{p}d_\Omega(x)^p$, and $\Vert x^* \Vert = \mu d_\Omega(x)^{p-1}$ for all $x^* \in \partialL  f(x)$.
%\begin{equation*}
%d_{{\rm{argmin}} f}(x) = d_\Omega(x), \quad
%f(x) - \inf f = \frac{\mu}{p}d_\Omega(x)^p,
%\quad \text{ and } \quad 
%\Vert x^* \Vert = \mu d_\Omega(x)^{p-1},
%\end{equation*}
%for all $x^* \in \partialL  f(x)$.
\end{proof}

\section{Distance bounds for globally Łojasiewicz functions}

In the previous section we have seen that square distance functions are examples of PŁ functions.
In this section we provide bounds justifying the idea that \emph{a PŁ function is a function growing like the square distance to its minimizers}.

\begin{theorem}\label{T:PL implies conditioning}
Let $f : \mathbb{R}^N \to \mathbb{R}$ be continuous and such that ${\rm{argmin}}~f \neq \emptyset$.
Let $p>1$, $\mu>0$, and let $f$ be globally $p$-Łojasiewicz with constant $\mu$.
Then $f$ is globally $p$-conditioned with constant $\mu/(p-1)^{p-1}$:
\begin{equation*}
(\forall x \in \mathbb{R}^N) \quad 
\frac{1}{(p-1)^{p-1}}\frac{\mu}{p} d_{{\rm{argmin}}f}(x)^p
\leq 
f(x) - \inf f.
\end{equation*}
\end{theorem}

\begin{proof}
This result comes from \cite[Theorem 24]{BolDanLeyMaz10}, where the authors need a few extra assumptions (which they do for simplifying their presentation). 
Since we drop those assumptions, we give here the details of the proof.
Let $x \in \mathbb{R}^N \setminus {\rm{argmin}}~f$ be fixed, and let us consider a sequence $(x_k)_{k \in \mathbb{N}}$ satisfying $x_0 = x$ and 
\begin{equation*}
x_{k+1} \in \underset{u \in \mathbb{R}^N}{\rm{argmin}}~f(u) + \frac{1}{2}\Vert u - x_k \Vert^2.
\end{equation*}
Note that our assumption that $f$ is bounded from below guarantees that $f + \frac{1}{2}\Vert \cdot - x_k \Vert^2$ is coercive.
Together with the continuity of $f$, we conclude that such sequence always exists.
Following \cite[Lemma 23]{BolDanLeyMaz10} and using the definition of our sequence, we write for every $k \in \mathbb{N}$ and every $u \in \mathbb{R}^N$ that
\begin{equation*}
\frac{1}{2}\Vert x_{k+1} - x_k \Vert^2 \leq \frac{1}{2}\Vert u - x_k \Vert^2 + f(u) - f(x_{k+1}).
\end{equation*}
Taking the infimum over the level set $S_{x_{k+1}} := \{ u \in \mathbb{R}^N \ | \ f(u) \leq f(x_{k+1}) \}$, and noting that $x_{k+1} \in S_{x_{k+1}}$, we deduce that
\begin{equation*}
\Vert x_{k+1} - x_k \Vert = d_{S_{x_{k+1}}}(x_k).
\end{equation*}
Now let $\varphi(t):= \frac{p}{q^{\frac{1}{q}} \mu^{\frac{1}{p}}} t^{\frac{1}{p}}$ with $\varphi'(t) = \frac{1}{q^{\frac{1}{q}} \mu^{\frac{1}{p}}} t^{\frac{-1}{q}}$ and $\varphi^{-1}(t) = \frac{1}{(p-1)^{p-1}}\frac{\mu}{p}t^p$.
It is easy to see that our global $p$-Łojasiewicz assumption on $f$ can be rewritten as
$1 \leq \varphi'(f(x) - \inf f) \Vert x^* \Vert$, for $x^* \in \partialL f(x)$.
Using a chain rule from \cite[Proposition 10.19.b]{RocWet09}, we know that 
\begin{equation*}
\partialL (\varphi \circ (f - \inf f))(x) \subset \varphi'(f(x) - \inf f) \partialL f(x).
\end{equation*}
We can use now \cite[Proposition 4.6]{DruIofLew15} to deduce, with the notations of \citep{BolDanLeyMaz10}, that $\vert \nabla (\varphi \circ (f - \inf f))\vert(x) \geq 1$, which means that the assumptions of \cite[Corollary 4]{BolDanLeyMaz10} are satisfied, and that we can now write
\begin{equation*}
\Vert x_{k+1} - x_k \Vert = d_{S_{x_{k+1}}}(x_k) \leq  \varphi(f(x_k) - \inf f) - \varphi(f(x_{k+1}) - \inf f).
\end{equation*}
Now use the fact that $x_0=x$, with the triangular inequality, to write
\begin{eqnarray*}
\Vert x_k - x \Vert 
& \leq & 
\sum\limits_{k=0}^{K-1} \Vert x_{k+1} - x_k \Vert 
\leq
\sum\limits_{k=0}^{K-1} \varphi(f(x_k) - \inf f) - \varphi(f(x_{k+1}) - \inf f) \\
& \leq & \varphi(f(x_0) - \inf f) - \varphi(f(x_{K}) - \inf f)
\leq \varphi(f(x) - \inf f).
\end{eqnarray*}
On the one hand, the above inequality implies that the sequence $x_k$ has finite length, which in particular means that it converges to some $x_\infty \in \mathbb{R}^N$, when $k \to + \infty$.
On the second hand, we can pass to the limit on $k$ to obtain
\begin{equation*}
\Vert x_\infty  - x \Vert \leq \varphi(f(x) - \inf f).
\end{equation*}
To conclude the proof, it remains to prove that $x_\infty \in {\rm{argmin}}~f$.
To see this, apply first-order optimality conditions to the definition of $x_{k+1}$ to write 
\begin{equation*}
x^*_{k+1} := x_k - x_{k+1} \in \partialL f(x_{k+1}),
\end{equation*}
and pass to the limit in the global $p$-Łojasiewicz inequality to obtain
\begin{equation*}
f(x_\infty) - \inf f 
= 
\lim\limits_{k \to + \infty} f(x_{k+1}) - \inf f
\leq
\lim\limits_{k \to + \infty} \frac{1}{q \mu^{\frac{q}{p}}}\Vert x^*_{k+1} \Vert^q = 0,
\end{equation*}
where the last equality comes from the fact that $\Vert x_{k+1}^* \Vert=\Vert x_k - x_{k+1} \Vert \to 0$.
This proves the claim after writing
\begin{equation*}
\varphi^{-1}(d_{{\rm{argmin}}f}(x))
\leq 
\varphi^{-1}(\Vert x - x_\infty \Vert)
\leq 
f(x) - \inf f.
\end{equation*}
\end{proof}

Next corollaries illustrate the idea that PŁ functions essentially behave like square distance functions.

\begin{corollary}\label{T:PL implies 2 conditioning}
Let $f:\mathbb{R}^N \to \mathbb{R}$ be continuous and such that ${\rm{argmin}}~f \neq \emptyset$. 
Let $f$ be Polyak-Łojasiewicz with constant $\mu>0$.
Then $f$ is globally $2$-conditioned with constant $\mu$ :
\begin{equation*}
(\forall x \in \mathbb{R}^N) \quad 
\frac{\mu}{2} d_{{\rm{argmin}}f}(x)^2
\leq 
f(x) - \inf f.
\end{equation*}
\end{corollary}

\begin{proof}
Apply Theorem \ref{T:PL implies conditioning} with $p=2$.
\end{proof}

\begin{corollary}\label{T:PL smooth conditioning sandwich}
Let $f:\mathbb{R}^N \to \mathbb{R}$ be differentiable and such that ${\rm{argmin}}~f \neq \emptyset$. 
Let $f$ be Polyak-Łojasiewicz with constant $\mu>0$, and $L$-smooth\footnote{We say that $f$ is $L$-smooth if $\nabla f$ is $L$-Lipschitz continuous.} with constant $L>0$.
Then $f$ behaves like a squared distance function, in the sense that
\begin{equation*}
(\forall x \in \mathbb{R}^N) \quad 
\frac{\mu}{2} d_{{\rm{argmin}}f}(x)^2
\leq 
f(x) - \inf f
\leq 
\frac{L}{2} d_{{\rm{argmin}}f}(x)^2
.
\end{equation*}
\end{corollary}

\begin{proof}
The first inequality was proved in the previous Corollary. For the second inequality, use the smoothness of $f$ (see e.g. Lemma 1.2.3 in \cite{Nes04}) to write for all $x^* \in {\rm{argmin}}~f$ that
\begin{equation*}
f(x) - f(x^*) - \langle \nabla f(x^*), x-x^* \rangle  \leq \frac{L}{2} \Vert x-x^* \Vert^2.
\end{equation*}
The conclusion follows from the fact that $f(x^*)=\inf f$, $\nabla f(x^*)=0$ and by taking $x^* = \proj_{{\rm{argmin}} f}(x)$.
\end{proof}

\section{Advanced comments}

\begin{remark}[On the continuity of $f$]
In Theorem \ref{T:PL implies conditioning} we assume that $f : \mathbb{R}^N \to \mathbb{R}$ is continuous. 
It is a simple exercise to show that the result is still true if we simply assume that $f : \mathbb{R}^N \to \mathbb{R} \cup \{+\infty\}$ is proper and lower semi-continuous.
\end{remark}

\begin{remark}[On the choice of the subdifferential]\label{R:choice of subdifferential}
In Definitions \ref{D:PL} and \ref{D:Łojasiewicz global} we require the Łojasiewicz inequalities to hold for every $x^* \in \partialL  f(x)$.
We could have considered a definition involving the Fréchet subdifferential $\partialF f(x) \subset \partialL  f(x)$ instead.
But the larger the subdifferential, the stronger the result, and we are indeed able to prove our main result in terms of the Limiting subdifferential.
One could wonder if the result remains true if we consider yet a larger subdifferential, such as the Clarke subdifferential \citep{Cla90}. 
As can be seen in the next counterexample, the answer is no, due to the fact that this Clarke subdifferential can detect sharp local maxima and contain $0$ at a point which is not a minimizer (a property incompatible with the PŁ definition). 
This is not surprizing since it is well-known that the Clarke subdifferential is sometimes \emph{too big} (see e.g. the discussion in Section 1.4.4 of \cite{Mor12}).
\end{remark}

\begin{figure}
\label{F:PL circle}
\begin{center}
\includegraphics[width=0.5\linewidth]{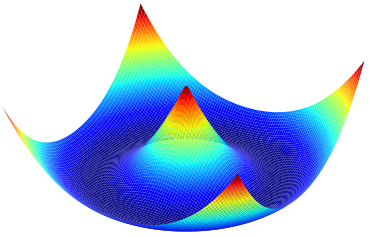}
\end{center}
\caption{Graph of $f(x) = \frac{\mu}{2} d_\Omega(x)^2$, where $\mu=20$ and $\Omega \subset \mathbb{R}^2$ is the unit circle. We can see a (nonsmooth) local maximizer at the origin.}
\end{figure}

\begin{example}[PŁ can fail  with Clarke subgradients]\label{Ex:PL clarke subgradients}
Let $\Omega = \mathbb{S}^2$ be the unit circle in $\mathbb{R}^2$, and $f(x) = \frac{1}{2}d_\Omega(x)^2$ (see Figure \ref{F:PL circle}).
At $x=0$, we  compute $f(x) - \inf f =\frac{1}{2}$, and (see Lemma \ref{L:subdifferential p distance}) $\partialL  f(x) = x - \proj_\Omega(x) = 0- \mathbb{S} = \mathbb{S}$.
Note how $x=0$ is not a Limiting critical point, even though it is a local maximum!
We recover as expected the fact that with $\mu=1$ we have
\begin{equation*}
f(0) - \inf f = \frac{1}{2}  \leq \frac{1}{2\mu} \Vert x^* \Vert^2
\text{ for all } x^* \in \partialL  f(0) = \mathbb{S}.
\end{equation*}
If we now compute the Clarke subdifferential at $x=0$, we obtain $\partialC  f(0) = \clco \partialL  f(0)$ (use Theorem 3.57 from \cite{Mor12} with the fact that $f$ is locally Lipschitz continuous).
That is, $\partialC  f(0)$ is the unit ball $\mathbb{B}$, which contains the origin.
This means that the inequality
\begin{equation*}
f(0) - \inf f = \frac{1}{2}  \leq \frac{1}{2\mu} \Vert x^* \Vert^2
\text{ for all } x^* \in \partialC  f(0) = \mathbb{B}
\end{equation*}
cannot be true when $x^* = 0$, whatever constant $\mu>0$ we consider.
The problem here is that $x=0$ is a \emph{Clarke} critical point which is not minimizer.
\end{example}

% Acknowledgements and Disclosure of Funding should go at the end, before appendices and references

\acks{The author would like to thank Robert Gower who accepted to proofread this document.}

\bibliography{references}

\end{document}